\documentclass{amsart}
\usepackage{amssymb}
\usepackage{upref}

\usepackage[all,cmtip]{xy}

\usepackage[lite,initials,shortalphabetic]{amsrefs}

\makeatletter
  \renewcommand{\p@enumi}{\thesubsection}
\makeatother

\newenvironment{resumeenumerate}[1]
{\begin{enumerate}
 \setcounter{enumi}{#1}
 \addtocounter{enumi}{-1}
}
{\end{enumerate}
}

\newenvironment{lettered}
{\begin{list}{\thelettercounter)}
 {\usecounter{lettercounter}\def\makelabel##1{\hss\llap{##1}}}
}
{\end{list}
}
\newcounter{lettercounter}
\renewcommand{\thelettercounter}{\alph{lettercounter}}

\makeatletter
\newcommand{\emsection}[1]{%
  \par
  \addpenalty\@secpenalty
  \vskip 6 pt plus 9 pt
  \emph{#1.}\nobreak\enspace\ignorespaces
}
\makeatother

\newcommand{\intro}{%
  \goodbreak
  \vskip 6 pt plus 9 pt
}

\numberwithin{equation}{subsection}

\newcommand{\Period}{\rlap{\enspace .}}

\newcommand{\cat}[1]{\boldsymbol{#1}}

\newcommand{\RelCat}{\mathbf{RelCat}}
\newcommand{\Cat}{\mathbf{Cat}}

\newcommand{\simp}{\mathrm{s}}

\newcommand{\spacedcdots}{{\cdot\;\cdot\;\cdot}}

\newcommand{\adjarrows}{\mathchoice{\longleftrightarrow}
  {\leftrightarrow}
  {\leftrightarrow}
  {\leftrightarrow}}

\begin{document}

\title[Quasi-categories and relative categories]
{A Thomason-like Quillen equivalence between quasi-categories and
  relative categories}

\author{C. Barwick}
\address{Department of Mathematics, Massachusetts Institute of
  Technology, Cambridge, MA 02139}
\email{clarkbar@math.mit.edu}

\author{D.M. Kan}
\address{Department of Mathematics, Massachusetts Institute of
  Technology, Cambridge, MA 02139}

\date{\today}

\begin{abstract}
  We describe a Quillen equivalence between \emph{quasi-categories}
  and \emph{relative categories} which is surprisingly similar to
  Thomason's Quillen equivalences between \emph{simplicial sets} and
  \emph{categories}.
\end{abstract}

\maketitle

%--------------------------------------------------------------------
%--------------------------------------------------------------------
\section{Introduction}
\label{sec:Intro}

In \cite{JT} Joyal and Tierney constructed a Quillen equivalence
\begin{displaymath}
  \cat S \adjarrows \simp\cat S
\end{displaymath}
between the Joyal structure on the category $\cat S$ of small
simplicial sets and the Rezk structure on the category $\simp\cat S$
of simplicial spaces (i.e.\ bi-simplicial sets) and in \cite{BK} we
described a Quillen equivalence
\begin{displaymath}
  \simp\cat S \adjarrows \RelCat
\end{displaymath}
between the Rezk structure on $\simp\cat S$ and the induced Rezk
structure on the category $\RelCat$ of relative categories.

In this note we observe that the resulting composite Quillen equivalence
\begin{displaymath}
  \cat S \adjarrows \RelCat
\end{displaymath}
admits a description which is almost identical to that of Thomason's
\cite{T} Quillen equivalence
\begin{displaymath}
  \cat S \adjarrows \Cat
\end{displaymath}
between the classical structure on $\cat S$ and the induced on on the
category $\Cat$ of small categories, as reformulated in
\cite{BK}*{6.7}.

\intro
To do this we recall from \cite{BK}*{4.2 and 4.5} the notion of
%--------------------------------------------------------------------
\section{The two-fold subdivision of a relative poset}
\label{sec:2foldsbdv}

For every $n \ge 0$, let
\begin{displaymath}
  \check{\cat n} \qquad\text{(resp. $\hat{\cat n}$)}
\end{displaymath}
denote the relative poset which has as underlying category the category
\begin{displaymath}
  0 \longrightarrow \spacedcdots \longrightarrow n
\end{displaymath}
and in which the weak equivalences are only the \emph{identity maps}
(resp.\ \emph{all maps}).

Given a relative poset $\cat P$, its \textbf{terminal} (resp.\
\textbf{initial}) \textbf{subdivision} then is the relative poset
$\xi_{t}\cat P$ (resp.\ $\xi_{i}\cat P$) which has
\begin{enumerate}
\item as \emph{objects} the \emph{monomorphisms}
  \begin{displaymath}
    \check{\cat n} \longrightarrow \cat P \qquad\text{($n \ge 0$)}
  \end{displaymath}
\item as \emph{maps}
  \begin{align*}
    (x_{1}\colon \check{\cat n}_{1} \to \cat P)&\longrightarrow
    (x_{2}\colon \check{\cat n}_{2} \to \cat P)\\
    \text{(resp. }(x_{2}\colon \check{\cat n}_{2} \to \cat P)
    &\longrightarrow (x_{1}\colon \check{\cat n}_{1} \to \cat P)\text{)}
  \end{align*}
  the commutative diagrams of the form
  \begin{displaymath}
    \xymatrix{
      {\check{\cat n}_{1}} \ar[rr] \ar[dr]_{x_{1}}
      && {\check{\cat n}_{2}} \ar[dl]^{x_{2}}\\
      & {\cat P}
    }
  \end{displaymath}
\end{enumerate}
and
\begin{resumeenumerate}{3}
  \item as \emph{weak equivalences} those of the above diagrams in
    which the induced map
    \begin{displaymath}
      x_{1}n_{1} \longrightarrow x_{2}n_{2}
      \qquad\text{(resp.\ $x_{2}0 \longrightarrow x_{1}0$)}
    \end{displaymath}
    is a weak equivalence in $\cat P$.
\end{resumeenumerate}

The \textbf{two-fold subdivision} of $\cat P$ then is the relative poset
\begin{displaymath}
  \xi\cat P = \xi_{t}\xi_{i}\cat P \Period
\end{displaymath}

%--------------------------------------------------------------------
%--------------------------------------------------------------------
\section{Conclusion}
\label{sec:conc}

In view of \cite{JT}*{4.1} and \cite{BK}*{5.2} we now can state:
\begin{enumerate}
\item \emph{the left adjoint in the above composite Quillen equivalence}
  \begin{displaymath}
    \cat S \adjarrows \RelCat
  \end{displaymath}
  \emph{is the colimit preserving functor which for every integer $n
    \ge 0$ sends}
  \begin{displaymath}
    \Delta[n] \in \cat S \text{\qquad to\qquad}
    \xi\check{\cat n} \in \RelCat
  \end{displaymath}
  \emph{and the right adjoint sends an object $X \in \RelCat$ to the
    simplicial set which in dimension $n$ ($n \ge 0$) consists of the
    maps $\xi\check{\cat n} \to X \in \RelCat$.}
\end{enumerate}
while, in view of the fact that $\Cat$ is canonically isomorphic to
the full subcategory
\begin{displaymath}
  \widehat\Cat \in\RelCat
\end{displaymath}
spanned by the relative categories in which \emph{every map is a weak
  equivalence} and \cite{BK}*{6.7},
\begin{resumeenumerate}{2}
\item \emph{the left adjoint in Thomason's Quillen equivalence}
  \begin{displaymath}
    \cat S \adjarrows \widehat\Cat
  \end{displaymath}
  \emph{is the colimit preserving functor which, for every integer $n \ge
  0$, sends}
  \begin{displaymath}
    \text{$\Delta[n]\in\cat S$ \emph{ to }
    $\xi\hat{\cat n} \in\widehat\Cat$}
    \qquad \text{($n \ge 0$)}
  \end{displaymath}
  \emph{while the right adjoint sends an object $X \in \widehat\Cat$
    to the simplicial set which in dimension $n$ ($n \ge 0$) consists
    of the maps}
  \begin{displaymath}
    \xi\hat{\cat n} \longrightarrow X \in \widehat\Cat
  \end{displaymath}
\end{resumeenumerate}

%--------------------------------------------------------------------
%--------------------------------------------------------------------

\end{document}